\documentclass[graybox]{svmult}

\usepackage{type1cm}        
\usepackage{makeidx}         
\usepackage{graphicx}        
\usepackage{multicol}        
\usepackage[bottom]{footmisc}

\usepackage{newtxtext}       %
\usepackage[varvw]{newtxmath}       

\makeindex             

\usepackage{amsmath}
\usepackage{amscd}
\usepackage{amsfonts}
\usepackage{dsfont}
\usepackage{cases}
\usepackage[applemac]{inputenc}
\usepackage{enumerate}
\usepackage[all]{xy}
\usepackage{esint}
\usepackage{url}

\newcommand{\R}{\mathbb R}

\newcommand{\F}{\cal F}
\newcommand{\G}{\Gamma}

\begin{document}

\title*{Interactions between resource dependent branching processes and equilibria}

 \titlerunning{Interacting branching processes} 
 
\author{F. Thomas Bruss}

\authorrunning{Interacting branching processes} 
\institute{F. Thomas Bruss \at Universit\'e libre de Bruxelles, Facult\'e des sciences, 
D\'epartement de Math\'ematique, CP 210, 
 Brussels, Belgium, \email{thomas.bruss@ulb.be}}
\maketitle

\abstract{This paper is a supplement to 
the paper "Interactions between Human Populations and Related Problems of Optimal Transport" written by the same author in honour of Marc Hallin, Universit\'e Libre de Bruxelles, at the occasion of Hallin's $75$th birthday. It was announced in the  main paper (Bruss (2024)) published in the Springer Festschrift entitled {\it Recent Advances in Econometrics and Statistics}. It contains the proofs which, given the space constraints required for the Festschrift, could not appear in the main paper. Moreover, we complement in the present supplement the main paper by brief comments on  related problems which are likely to turn up in practice  for problems of guiding human populations, namely problems of control and problems of optimal stopping.}

\smallskip

\noindent Math. Subj. Class.:  60J85
\subsection{Particularities of resource dependent branching processes}
As far as we know, the first version of a resource dependent branching process was introduced by the author in his talk
at Cornell University during the Conference on Stochastic Processes and their Applications  (1983). It was a model
of a branching process where in each generation particles reproduced according to a Galton-Watson process, provided that they were 
given sufficient resources by the preceding generation. 

In this model, particles used resources within their generation, and they left resources and/or created new resources for the next generation. A  rule established by a {\it society} decided how to distribute resources. Particles which did not receive the quantity of resources according to their needs were supposed to emigrate without leaving offspring. The focus was first on extinction or survival only. Criteria for possible survival of the population could still profit from those found in Bruss (1978) generalising a model studied by Zubkov (1970). At the beginning, only a single population was considered, Hence there was not yet a question concerning possible equilibria between sub-populations. The same was still true for the more elaborate model of
Bruss and Duerinckx (2015). Now the challenge is to find tractable models for populations which consist of
interacting sub-populations.

When considering  a population consisting of sub-populations, then equilibrium questions turn up in a natural way.  Tractability of such models seem to depend strongly on the assumption that within each sub-population, reproduction is again according to the scheme of a Galton-Watson process. To understand the main reasons why,
it may be helpful to refer to the paper of  Bruss (2022), which underlines the particular role of Galton-Watson processes as building blocks for more complicated branching processes. 

\subsection{Bisexual reproduction}
Before going to the main proofs, we must recall that
the Bruss-Duerinckx model (2015) of resource dependent branching processes is somewhat simplified as far as the bisexual reproduction of (human) individuals is concerned. Namely the reproduction  
rate is seen as the {\it average reproduction rate of mating units} as defined in Bruss (1984). This simplification is justified for all objectives and results in Bruss and Duerinckx (2015) as well as in the present supplement.

It is important however  to note that, more generally, the interaction of mating functions affects the branching property for the (true) individuals as such. Consequently,  for certain questions concerning bisexual population processes the introduction of {\it averaged reproduction means} of mating units  may not be adequate,
as in Daley (1968) but also  in Alsmeyer and R\"osler, Molina et al. (2002),  Molina et al. (2007) and others. If the notion of the {\it average reproduction rate of mating units} is not suitable but the interest of the study stays focussed on a limiting behaviour then we refer to Bansaye et al. (2023)
who found interesting new results on what these authors call the relevant {\it scaling limits} for bisexual branching processes. 

Results on multi-type processes are given in Hautphenne et al.(2013) and Braunsteins et al. (2022), and a new study of  bisexual processes can be found in Fritsch et al. (2024). 

\subsection{Asymptotic equilibrium between two sub-populations}

Here is now the central result for a society consisting of two sub-populations. 

We exemplify our two populations by what we call {\it  home population} and {\it immigrant population,} from which individuals migrate into the home-population.
We recall that the parameter $m$ refers to the mean reproduction rate or natality rate of individuals within a given population whereas $r$ denotes the mean resource production of a randomly chosen individual in the same population. $F$ denotes the claim-size distribution function of individuals within the given population, respectively sub-population. The indices $h$ respectively $i$ are mnemonic for {\it home} respectively {\it 
immigrant.} The two claim-size distribution functions
$F_h$ and $F_i$ are supposed to be absolute continuous.
\subsection{Proof of Theorem 3 ({\rm Bruss (2024))}} 
The following is the complete proof of Theorem 3 of Bruss (2024).  
In the latter, the indices $h$ and $i$ refer to the home-population and to the immigrant population, respectively, that is $F_h$ ($F_i$) denotes the claim distribution function of the home population (immigrant population), $r_h$ ($r_i$) and $m_h$ ($m_i$) the resource production mean and the reproduction mean of a randomly chosen individual
within these sub-populations, respectively. The effectives in each sub-population at time $t$ are correspondingly denoted by $\G^h_t$ and $\G^i_t.$
 
 The proof consists of four parts, the first three  (I)-(III) proving (a), and (IV) proving (b) of Theorem 3.

\medskip

\noindent (I)~We will  first confirm our intuition that, if a limiting equilibrium exists, then necessarily both sub-processes tend to infinity as $t\to \infty$, i.e.\begin{align} \label{ratio}\lim _{t \to \infty}\left(\G^i_t/\G^h_t \right){\rm exists~}\implies P(\G^i_t \to \infty)|\G^i_t \not \to 0)\,=\,P(\G^h_t \to \infty| \G^h_t \not \to 0)=1.\end{align}

\smallskip 
 
The proof of part (I) is by contradiction.

Suppose that the  statement (\ref{ratio}) is wrong. Clearly both processes $(\G^i_t)$ and $(\G^h_t)$ must stay bounded away from zero as $t\to\infty,$ because the state zero is an absorbing state, given that there are no new immigrants after time $0.$ Furthermore, for $\alpha$ to be meaningful as an equilibrium between the two sub-populations, $\alpha$ must satisfy $0<\alpha<\infty$, that is, neither of the sub-populations $(\G^i_t)$ nor  $(\G^h_t)$ may suffer extinction or be "wipe out" by the other one. 

It follows that if our statement were wrong, there would exist bounds $b_h>0$ and $b_i>0$, say, such that $$P(\G^i_t \le b_i ~\rm{i.o.}| \G^i_t > 0) >0 ~\rm {and/or~}P(\G^h_t \le b_h ~\rm{i.o.}| \G^h_t > 0) >0,$$ where i.o. is shorthand for {\it infinitely often.} Now put $b=\max\{b_i, b_h\}$ and $p_0=\min\{p_0^i, p_0^h\},$  where $p^h_0>0$, respectively $p^i_0>0$, denotes the probability that a randomly chosen individual in the home-population, respectively immigrant-population, has no offspring. Since reproduction  is mutually independent for all individuals within the same sub-population, we must have \begin{align}\sum_t
P(\G^h_{t+1}=0|\G^h_t>0) =\infty {~~\rm or~~}\sum_t
P(\G^i_{t+1}=0|\G^i_t>0) = \infty\end{align}because all terms are non-negative in both sums, and, in at least one sum, infinitely many terms are not less than $p_0^b>0. $ But the latter implies  (see e.g. Lemma 1 and Corollary 1 of Bruss (1980)) that at least one sub-process will die out almost surely. This contradicts the assumption $0<\alpha<\infty$, and hence (\ref{ratio}) is proved.

\smallskip
\noindent (II)~We now prove that, for a given $0<\alpha<\infty$ satisfying the definition of an equilibrium, there must exist a value $\tau:=\tau(\alpha)$ such that equation (14) in the main paper (Bruss (2014)) is satisfied.
We first note that if such a value $\tau$ exists for a given value $\alpha,$ then $\tau$ is unique if  the densities $dF_h(t)/dt$ and $dF_i(t)/dt$ do not vanish at the same time in a neighbourhood of $\tau.$ Indeed, this follows since both integrands on the l.h.s. of the equilibrium equation (see (14) in the main paper))  are non-negative. If we now denote $\cal S$ the union of the supports of $F_h$ and $F_i$ we can define more generally
\begin{align}\tau:=\inf\left\{t\in {\cal S}: m_h\int_0^t xdF_h(x) + \alpha\,m_i\int_0^t xdF_i(x) =r_h +\alpha r_i\right\}.\end{align} This implies the uniqueness of $\tau$ in any case, and justifies the notation $\tau:=\tau(\alpha).$

We now turn to equation (12) in the main paper under the assumption  $\tilde R(\G_t^h,  \G^i_t) = R^h_t(\G_t^h)+R^i_t(\G^i_t)$ which follows from our hypothesis that, as we recall, all resource claims are submitted to the same accumulated resource space.
Also, since reproduction rates and resource production rates of individuals are independent  within each sub-process, and  $F_h$ and $F_i$ are fixed distribution functions, we can apply the strong law of large numbers in equation (12) via $\G^h_t$ for both processes separately. At the same time we will see that $\tau_t$ converges to a constant $\tau:=\tau(\alpha)$ almost surely.

Indeed, by dividing both sides of equation (12) in the main paper by $\G^h_t,$  and using the dummy multiplication factor $\G_t^i/\G_t^i$ for the second terms on both sides, we obtain\begin{align} 
\frac{D_t^h(\G^h_t)}{\G^h_t}\int_0^{\tau_t}xdF_h(x)+
\frac{D_t^i(\G^i_t)}{\G^i_t}\frac{\G^i_t}{\G^h_t}\int_0^{\tau_t}xdF_i(x)=\frac{R_t^h(\G_t^h)}{ \G_t^h}+\frac{R^i_t(\G^i_t)}{\G^i_t}\frac{\G^i_t}{\G_t^h}.\end{align}
Conditioned on survival of both sub-processes $(\G_t^h)$ and $(\G_t^i)$ we know from the proof of part (I) and the strong law of large numbers that the term multiplying the first integral on the l.h.s. converges almost surely to $m_h$  whereas the first term on the r.h.s. converges almost surely to $r_h.$  Moreover, if the limit $\alpha$ exists then, as $t \to \infty,$ we must have
\begin{align}\frac{D_t^i(\G^i_t)}{\G^i_t}\frac{\G^i_t}{\G^h_t}\to \alpha\, m_i {\rm~~ a.s.}~{\rm ~and}~\frac{R^i_t(\G^i_t)}{\G^i_t}\frac{\G^i_t}{\G_t^h}\to \alpha\, r_i {\rm~~ a.s.}\end{align}
This implies that, conditioned on survival of both sub-populations and on the existence of the limit $\alpha$, the r.h.s. of the obtained equation has the limit $r_h+\alpha r_i$ a.s. as $t \to \infty.$ Consequently, its corresponding l.h.s. must  also have  a  limit.  Since the upper bound in both integrals is the same value $\tau_t$ in (4) above, and both integrals have non-negative integrands,  $\tau_t$ must converge almost surely  to a constant $\tau$ as $t\to \infty$. 

Taking now the preceding two arguments together we conclude that, if an equilibrium exists, then the  corresponding $\alpha$ and $\tau$ must satisfy the limiting analogue of the limiting equation in (5) above, namely 
$$m_h\int_0^\tau xdF_h(x) + \alpha\,m_i\int_0^\tau xdF_i(x) =r_h +\alpha r_i.$$ This is equation (14)  as claimed in the Theorem 3.
Moreover, with our definition of $\tau$ for a given $\alpha$, the $\tau:=\tau(\alpha)$  must  be unique. 
This proves part (II).

\medskip\noindent
(III) We now have to prove that the combined constraint qualifications in part (a), namely $m_h\, F_h(\tau)=m_i\, F_i(\tau) \ge1,$  must hold for $\tau=\tau(\alpha).$ For this to see we first prove the  equality part of it, which, at a first look, may seem very restrictive. Let the random variable $\alpha_t$  be defined by $\alpha_t=\G^i_t/\G^h_t,$ or equivalently\begin{align} \frac{1}{1+\alpha_t}=\frac{\G_t^h}{\G_t^h+\G_t^i}, ~ t=1,2 \cdots.\end{align}
Conditioned on survival of both sub-processes we  know thus from part (I) that $\alpha_t\to\alpha$ for some $\alpha$ and, as seen in part (II), and $\tau_t \to \tau:=\tau(\alpha)$ a.s. as $t\to \infty.$ 
Further, in order to be an element of the home-population at time $t+1,$ it is necessary for an individual to be a descendant of it, and if so then it is sufficient if its resource claim does not exceed the threshold $\tau_{t}.$ 
It follows that, conditioned on survival of  $(\G_t^h)$ and $(\G_t^i),$ the random fraction of individuals belonging to
the home-population one generation later can  be written as
\begin{align}\frac{1}{1+\alpha_{t+1}}=\frac{D_t^h(\G_t^h) F_h(\tau_t)}{D_t^h(\G_t^h) F_h(\tau_t)+D_t^i(\G_t^i) F_i(\tau_t)},  \end{align}
where $\tau_t \to \tau$ a.s. as $t \to \infty.$
We now divide on the r.h.s. of this equation the numerator and denominator by $\G_t^h.$ (Recall that the latter must be positive.) Using part (I) of the proof and the existence of the limit $\alpha$ we see  then that, conditioned on survival of both sub-processes,  $$D^i_t(\G_t^i)/D^h_t(\G_t^h)\to m_i \alpha /m_h  \rm{~a. s. ~as~} t \to \infty.$$ Taking the limit on both sides of (8) above for $t\to \infty$ yields then after straightforward computations
\begin{align} \frac{1}{1+\alpha}=\frac{1}{1+\alpha \left(m_i F_i(\tau)/m_hF_h(\tau)\right)},\end{align}
and hence $m_i F_i(\tau)=m_hF_h(\tau).$ 
This proves the equality part of the constraint qualification (15).

\medskip To complete the proof of part (III) it remains to show that  the conditions $m_hF_h(\tau) \ge 1$ and $m_hF_h(\tau) \ge 1$  are necessary for the existence of an equilibrium.
We argue again by contradiction. 

\smallskip\noindent
Suppose the contrary, and let us first suppose that $0\le m_hF_h(\tau) < 1.$ Let $b \in \,] m_hF_h(\tau), 1[.$  Since  we can confine our interest on the case $\G_t^h \to \infty$ almost surely as $t \to \infty$, and since $D_t^h(\G_t^h)$ is the sum of $\G_t^h$ i.i.d. random variables, we see straightforwardly from the strong law of large numbers and $b<1$ that, for all $\G_t^h$ sufficiently large, \begin{align}E\left(\G_{t+1}^h\big|\G_t^h>0\right) = E\left(D_t^h(\G_t^h)F_h(\tau_t)\big|\G_t^h>0\right)\\< E\left(b\G_t^h\big| \G_t^h>0\right)<E\left(\G_t^h\big| \G_t^h>0\right)\nonumber.\end{align} This implies that the process $(\G_t^h)$, conditioned on non-extinction, stays bounded in expectation so that \begin{align}\sum_t (p_o^h) ^{E(\G_{t+1}^h|\G_t>0)} = \infty,\end{align}where we recall that $p_o^h$ denotes the probability that an individual in the home-population has no children. But then, by another Borel-Cantelli type argument  (see e.g. Bruss (1978), pp 54-56, Theorem 1) we get the contradiction $(\G_t^h) \to 0 {\rm ~a.s.}$ as $t \to \infty.$

A contradiction is obtained in an analogous way by supposing that  an equilibrium exists and $m_iF_i(\tau) < 1.$ 
 This completes the proof of part (III).

\medskip \noindent
(IV) We finally have to show that, given that $\{\G_t^i/\G_t^h\not\to 0\}\cap \{\G_t^i/\G_t^h \not\to \infty\},$ then, replacing the constraint qualification  $m_h\, F_h(\tau)=m_i\, F_i(\tau) \ge 1$ by the slightly stronger condition \begin{align*}m_h\, F_h(\tau)=m_i\, F_i(\tau) >1\end{align*}shows that this is  sufficient  for an $\alpha$-equilibrium to exist with a strictly positive probability.  

Since we know that the equality
$m_h\, F_h(\tau)=m_i\, F_i(\tau)$  is necessary for the existence of an equilibrium, as we have shown already in the first part of the proof of part (III), we only have to show now that \begin{align}P\Big(\{\G_t^h \to \infty\}\cap\{\G_t^i\to \infty\}\Big|m_hF_h(\tau)=m_iF_i(\tau)>1\Big)>0,\end{align}
because, given the joint event $\{\G_t^h \to \infty\}\cap\{\G_t^i\to \infty\},$  we can follow the arguments from (18) up to (21) in the main paper to establish its limiting equation (15).

Now,  if at least one of the two sub-processes tends to infinity with strictly positive probability, then {\it both}
must do so according to the definition of $\alpha$ ($0<\alpha<\infty$)
as being the a.s. limiting ratio of $\G_t^i/\G_t^h$ as $t \to \infty.$ Hence, recalling  part (I), it suffices to show that 
\begin{align}P\Big(\G_t^h \to \infty \,\Big|\,m_hF_h(\tau)>1\Big)=1-P\Big(\G_t^h \to 0\,\Big|\,m_hF_h(\tau)>1\Big)>0.\end{align} Since, the process $(\G_t^h)$ is a RDBP  by definition, the latter follows from
Theorem 4.4 ii) b) of Bruss and Duerinckx  (2015).

This completes the proof of Part (VI) of the Theorem,
and thus all parts of the proof are complete.

\subsection{Several sub-populations} The paper Bruss (2024) does not treat the case of more than two sub-populations. However, using the so-called BRS-inequality (see Steele (2016) and Bruss (2021)) necessary survival and equilibrium results can be derived in an almost analogous way for an arbitrary number of sub-populations. As in the preceding proof, it needs the assumption that the reproduction of individuals stays independent within each sub-population.

Again, in order to treat equilibrium questions, one typically has to look first at extinction probabilities of processes,  questions which (under real world assumptions) are clearly linked with questions of how quickly the considered population or sub-population will grow. Hence relationships between the size of a population and its extinction probability attract interest. The more general they are, the more interest they typically attract. See e.g. the new paper by Ball et al. (2024)
on the general relationship between population size and risk of extinction.

Of course, in reality it is often the case that  certain sub-populations  integrate into another one.  For example, the sub-population of immigrants (or a part of it) often integrates into the home-population. We just mention here that the equilibrium problem
can then be solved in a similar way as before, if the newly integrated individuals behave exactly as the individuals
of the "host"-population, and if there are no new immigrants after time $0.$
\subsection{Optimal transport and the equilibrium equation for two sub-populations}

~~~~~On the search for reaching an equilibrium between sub-populations, decision makers are likely to have to look for control tools which are intended to affect certain distributions and/or certain parameters.  At the same time decision makers must usually confine their interest to those parts of the economic interplay, where the tools look to be most promising. In the setting described in our main paper, these are the claim distribution functions $F_h$ and $F_i$ as well as the expected resource production rates $r_h$ and $r_i.$ However, additional constraints may intervene on both the claim distribution functions and on the resource productivity (capacity) of individuals in the respective sub-populations. 

In the following we discuss one of the simpler cases where $r_h$ and $r_i$ are seen as fixed and where the control must be based on the claim distributions only.
But then, what does it mean to concentrate one's interest on the most promising part of the interplay? 

We can argue that our interest should first concentrate on the claim distribution of the "richer" sub-population,
that is, on the one allowing before control for the higher standard of living. Let us assume for instance that this one is $F_h$. Adapting the terminology often used in problems of optimal transport we call then $F_h$ the "supply" distribution function whereas we see the distribution function $\tilde F_h$ into which we would like to transport it as {\it a} "demand" distribution function. Since in practice we may have more than one choice we define $\tilde F_h$ to be any distribution function $F$ which allows for a solution $(\tilde\tau, \tilde\alpha)$ of 
 \begin{align} m_h\int_0^{\tilde \tau}x dF(x) +\tilde\alpha m_i\int_0^{\tilde\tau}x dF_i(x)dx=r_h+\tilde\alpha r_i,\end{align} subject to $m_hF(\tilde\tau)=m_iF_i(\tilde\tau)>1.$ (Compare with equation (14) in Theorem 3 of Bruss (2014).)
 
 Once a cost function for shifting one unit of mass from $x$ to $y$, denoted by $c(x,y),$ is determined
 we then have the optimal transport problem as presented by Rachev and R\"uschendorf (1998) in their chapter 7.1, which is as follows:
 
 Let for discrete distribution functions$F_h$ and $\tilde F_h$ the set of distribution functions $F$ on $\R^2$ with marginals $F_h$ and $\tilde F_h$ be denoted by ${\F }(F_h, \tilde F_h ).$ Our assumptions imply then $F(x,\infty)=F_h(x)$ and $F(\infty,y)=\tilde F(y).$ The control problem seen as an optimal transport problem with cost function $c(x,y)\ge 0$ for $x,y \in \R$ becomes
 \begin{align}\min_{F\in{\cal F}}~ \int_{\R^2} c(x,y) \,dF(x,y). \end{align}
 The preceding problem is the continuous version of the  discrete analogue problem
 \begin{align}{\rm Minimize } \sum_{i=1}^m\sum_{j=1}^n ~ c_{i,j}x_{i,j}  {\rm ~subject~ to}
 \end{align}\begin{align}\forall~ i,j:~\sum_{j=1}^n  x_{i,j}=a_i, 1\le j\le m; \sum_{i=1}^m x_{i j}=b_j,~ j=1,2,\cdots,n ~{\rm with~}x_{i j} \ge 0. \end{align}
If we suppose that the cost function $c{x,y}$ respectively the $c_{i j}$ satisfy the so-called Monge condition  then the solution is given by 
\begin{align} F^*(x,y) = min ~ \Big\{F_h(x), ~\tilde F_h(y)\Big\}
\end{align} 
The book by Rachev and R\"uschendorf (1998) treats several aspects of this problem in sufficient detail. See also Villani (2009) for a newer review of related results.

A large range of computational aspects of optimal transport problems is studied in  Peyr\'e and Cuturi (2019).)
\subsection{Aspects of control}
As far as the author understands, the discrete analog (15) and (16)
given above should be in practice more important for the intrinsic
control problem because it is not likely (and probably not even desirable) that decision makers would try to formulate the demand distribution function in "infinitesimal" precision. The reason is that decision makers are typically bound to think in terms of small changes
for which they can observe a visible impact. 

Now, if changes take time and cause costs, as they usually do, then certain envisaged changes may have to be given up at a  premature state. Keeping this in mind we see that some of the control problems may turn also into problems of Optimal Stopping. If so, then they are typically complicated and one may have to look for good approximations of the optimal solution. In this task of finding good approximations for the optimal stopping problem, several results of R\"uschendorf (2016)) should be rewarding.

 In our case, stopping means to give up on a certain control in order to replace a currently pursued objective by a
 modified one, or even a completely new one. Having said this, it does not seem to be easy to formulate the control problem in such a way that it combines the (final) goal of optimal transport and an intermediate (goal) of optimal stopping in a way that 
would always makes sense. 
Combining measure transportation and statistical decision theory, as studied in Hallin (2022), may yield an alternative fruitful approach to control resource dependent branching processes, but the author has no sufficient experience to judge. 
\subsection{More than two sub-populations} As far as questions of optimal transport
are concerned, problems become in general harder for more than two sub-populations. The difficulty becomes obvious if one recalls that changing one claim-distribution function (for a given sub-population) may then  have
an influence of several  other (possibly all) claim distribution functions. 

Hence, if there are $s$ sub-populations, one would have to consider
in the worst case (i.e. if a priori no case can be eliminated as being unreasonable) at an intimidating number of $2^s-1$ different models
and problems.
 \subsection*{~~~~References}
 ~~~~Alsmeyer  G. and R\"osler U. (2005)  {\it Asexual Versus Promiscuous Bisexual Galton-Watson Processes: The Extinction Probability Ratio,}
Ann. Appl. Probab., Vol. 12, No. 1, 125-142.

\smallskip 
Ball T. S.,  Balmford B.,  Balmford A., Rinaldo D., Visconti P. and Green R. (2024),
{\it A general relationship between population size and extinction risk}, arXiv:2411.13228v1.

\smallskip
Bansaye V., Caballero M.-E., M\'el\'eard S. and 
San Martin J. \,(2023), {\it Scaling limits of bisexual Galton-Watson processes}, Stochastics
Vol. 95, Issue 5, 749-784.

\smallskip
Braunsteins P., Decrouez G., and S. Hautphenne (2022). 
{\it A pathwise approach to the extinction of branching processes with countably many types},
Stoch.\ Proc.\ Appl., {\bf 129}, 713--739.

\smallskip
Bruss F.T. (1978), {\it Branching processes with random absorbing processes}, J. Appl. Prob. Vol. 15,
54-64.

\smallskip
Bruss F.T.   (1980), {\it A Counterpart of the Borel-Cantelli Lemma}, J. Appl. Prob., Vol. 17, 1094-1101.

\smallskip
Bruss F.T.   (1984) ~{\it A Note on Extinction Criteria for Bisexual Galton-Watson Processes}, J. Appl. Prob., Vol. 21, 915-919.

\smallskip
Bruss F.T. (2021) {\it The BRS-inequality and its applications},  Probab. Surveys, Vol. 18, 44-76.

\smallskip
Bruss F.T. (2022)~{\it Galton-Watson processes and their role as building blocks for branching processes}, Theory Probab.  Appl., Vol. 67, No 1, 141-153.

\smallskip
Bruss F.T. (2024), {Interactions between Human Populations and Related Problems of Optimal Transport,} in {\it Recent Advances in Econometrics and Statistics}, Springer Verlag.

\smallskip
Bruss F.T. and Duerinckx M. (2015) ~{\it Resource dependent branching processes and the envelope of societies},
 Ann. of Appl. Probab., Vol. 25, Nr 1, 324-372.
 
\smallskip 
Daley D.J. (1968){ \it Extinction conditions for certain bisexual Galton-Watson branching processes}, Z. Wahrscheinlichkeitsth. und Verw. Geb.,Vol. 9, Issue 4,  315-322. 

\smallskip
Fritsch C., Villemonais D., Zalduendo N. (2024) 
{\it The Multi-type Bisexual Galton-Watson Branching Process,}  Ann. Inst. H. Poincar\'e,  Probab. Statist. 60(4): 2975-3008.

\smallskip
 Hallin M. (2022)~{\it Measure Transportation and Statistical Decision Theory}, Annual Review of  Statistics and its Applic., Vol. 9: 401-424.
 
\smallskip
Hautphenne S., Latouche G. and  Nguyen G. (2013),
{\it Extinction probabilities of branching processes with countably infinitely many types},
Adv. in Appl. Probab. 45(4): 1068-1082. 

\smallskip Molina M., Mota M. and  Ramos A. (2002),{\it Bisexual Galton-Watson Branching Process with Population-Size-Dependent Mating}, J. Appl. Prob. Vol. 39, No. 3, 479-490.

\smallskip
Molina M., del Puerto I. and Ramos A. (2007), {\it A class of controlled bisexual branching processes with mating depending on the number of progenitor couples,} Stat. and Prob. Letters,
Vol. 77, Issue 18, 1737-1743.

\smallskip
Peyr\'e G. and Cuturi M. (2019), {\it Computational Optimal Transport: With Applications to Data Science}, In Foundations and Trends in Machine Learning: Vol. 11: No. 5-6, 355-607. 

\smallskip
Rachev S. and R\"uschendorf L. (1998). {\it Mass Transportation Problems}, 
Vol. II: Application, Springer-Verlag

\smallskip
R\"uschendorf \, L.\, (2016), \,{\it Approximative  solutions of optimal stopping and selection problems}, Math. Applicanda, Vol. 44(1), Tom 22/60,  17- 44.

\smallskip
Steele J.M. (2016), {\it The Bruss-Robertson Inequality: Elaborations, Extensions, and Applications}, Math. Applicanda, Vol. 44(1), Tom 22/60, 3-16.

\smallskip
Villani C. (2009), {\it Optimal Transport, Old and New}, Vol. 334 of Grundlehren der
Mathematischen Wissenschaften, Springer, New York.

\smallskip
Zubkov A.M. (1970), {\it A degeneracy condition for a bounded branching process}, Mat. Zametki
Vol. 8, 9-18; English translation: Math. Notes Vol. 8, 472-477,

 \end{document}